\documentclass[12pt,a4paper,reqno]{amsart}
\usepackage{amsmath,latexsym,epsfig,graphics,amsfonts,amssymb}

\newtheorem{theorem}{Theorem}[section]
\newtheorem{lemma}[theorem]{Lemma}
\newtheorem{proposition}[theorem]{Proposition}

\numberwithin{equation}{section}

\newcommand{\Z}{\mathbb{Z}}
\newcommand{\R}{\mathbb{R}}
\newcommand{\ThF}{F}
\newcommand{\ThT}{\textup{T}}
\newcommand{\ThV}{\textup{V}}

\begin{document}
\pagestyle{plain} \pagenumbering{arabic}
\title{Conjugacy, roots, and centralizers in Thompson's group $\ThF$}

\author{Nick Gill}
\email{nickgill@cantab.net}
\address{Institute of Mathematical Sciences, CIT Campus, Taramani, Chennai 600113, India}

\author{Ian Short}
\email{ian.short@nuim.ie}
\address{Mathematics Department, Logic House, N.U.I. Maynooth, Maynooth, County Kildare, Ireland}

\subjclass[2000]{Primary: 20E45. Secondary: 20F10, 37E05}

\keywords{Thompson, piecewise linear, conjugacy, roots, centralizers}

\date{\today}

\begin{abstract}
We complete the program begun by Brin and Squier of characterising conjugacy in Thompson's group  $\ThF$ using the standard action of $\ThF$ as a group of piecewise linear homeomorphisms of the unit interval. 
\end{abstract}

\maketitle
\section{Introduction}


The object of this paper is to extend the methods of Brin and Squier described in \cite{BrSq01} to solve the conjugacy problem in Thompson's group $\ThF$, and to analyse roots and centralizers in $\ThF$. 

Let $\text{PL}^+(a,b)$ denote the group of piecewise linear order-preserving homeomorphisms of an open interval $(a,b)$. The points at which an element $f$ of $\text{PL}^+(a,b)$ is not locally affine are called the \emph{nodes} of $f$. Thompson's group $\ThF$ is the subgroup of  $\text{PL}^+(0,1)$ defined as follows: an element $f$ of $\text{PL}^+(0,1)$ lies in $\ThF$ if and only if the nodes of $f$ lie in the ring of dyadic rational numbers, $\mathbb{Z}\left[\frac12\right]$, and $f'(x)$ is a power of $2$ whenever $x$ is not a node. 

In \cite{BrSq01} Brin and Squier analysed conjugacy in $\text{PL}^+(a,b)$ for $(a,b)$ any open interval. For $(a,b)$ equal to $(0,1)$ we can restate their primary result \cite[Theorem 5.3]{BrSq01} as follows: we have a simple quantity $\Sigma$ on $\text{PL}^+(0,1)$ such that two elements $f$ and $g$ of $\text{PL}^+(0,1)$ are conjugate if and only if $\Sigma_f=\Sigma_g$. If $f$ and $g$ are elements of $\ThF$ then $\Sigma_f$ and $\Sigma_g$ can be computed and compared using a simple algorithm. Brin and Squier comment on their construction of $\Sigma$ that,``Our goal at the time was to analyze the conjugacy problem in Thompson's group $\ThF$." In this paper we achieve Brin and Squier's goal by defining a quantity $\Delta$ on $\ThF$ such that the following theorem holds.

\begin{theorem}\label{T: main}
Let $f,g\in \ThF$. Then $f$ and $g$ are conjugate in $\ThF$ if and only if 
$$(\Sigma_f, \Delta_f)=(\Sigma_g, \Delta_g).$$
\end{theorem}

This is not the first solution of the conjugacy problem in $F$. In particular the conjugacy problem in $\ThF$ was first solved by Guba and Sapir in \cite{Gu97} using diagram groups.  More recently, Belk and Matucci \cite{belkmatt1,belkmatt2} have another solution using strand diagrams. Kassabov and Matucci \cite{kassmatt} also solved the conjugacy problem, the simultaneous conjugacy problem, and analysed roots and centralizers in $\ThF$. Our analysis is different to all of these as we build on the geometric invariants introduced by Brin and Squier.




We will not prove Theorem \ref{T: main} directly. Rather we prove the following proposition which, given \cite[Theorem 5.3]{BrSq01}, implies Theorem \ref{T: main}:

\begin{proposition}\label{P: main}
Two elements $f$ and $g$ of $F$ that are conjugate in $\textup{PL}^+(0,1)$ are conjugate in $F$ if and only if $\Delta_f=\Delta_g$.
\end{proposition}

Our paper is structured as follows. In \S\ref{S: PL} we introduce some important background concepts, including the definition of $\Sigma$. In \S\ref{S: Delta} we define $\Delta$. In \S\ref{S: invariant} we prove Proposition~\ref{P: main}. In \S\ref{S:calc1} and \S\ref{S:calc2} we outline formulae which can be used to calculate $\Delta$. Finally, in \S \ref{S:rc}, we discuss roots and centralizers in $\ThF$.

It is likely that the results in this paper can be extended to hold in Thompson's groups $\ThT$ and $\ThV$; we hope to describe such extensions in a forthcoming paper.

\section{The definition of $\Sigma$}\label{S: PL}

Let $f$ be a member of Thompson's group $\ThF$, embedded in $\text{PL}^+(0,1)$. Following Brin and Squier \cite{BrSq01} we define the invariant $\Sigma_f$ to be a tuple of three quantities, $\Sigma_1, \Sigma_2$ and $\Sigma_3$, which depend on $f$.

The first quantity, $\Sigma_1$, is a list of integers relating to values of the \emph{signature of} $f$, $\epsilon_f$. We define $\epsilon_f$ as follows:
$$\epsilon_f:\R\to\{-1,0,1\}, \ \ x\mapsto \left\{
\begin{array}{ll}
1, & f(x)>x; \\
0, &  f(x)=x; \\
-1, & f(x)<x.
\end{array}
\right.$$
If $f$ is an element of $\text{PL}^+(0,1)$ then there is a sequence of open intervals 
\begin{equation}\label{E: intervals}
I_1,I_2,\dots,I_m,\quad I_j=(p_{j-1},p_j),\quad p_0=0, \quad p_m=1,
\end{equation}
 such that $\epsilon_f$ is constant on each interval, and the values of $\epsilon_f$ on two consecutive intervals differ. We define $\Sigma_1=(\epsilon_f(x_1), \dots, \epsilon_f(x_m))$ where $x_i\in I_i$ for $i=1,\dots, m$.

Let $\text{fix}(f)$ be the set of fixed points of $f$ and observe that the points $p_0,\dots,p_m$ from (\ref{E: intervals}) all lie in $\text{fix}(f)$. We say that the interval $I_j$ is a \emph{bump domain} of $f$ if $\epsilon_f$ is non-zero on that interval. Our next two invariants consist of lists with entries for each bump domain of $f$.

If $k$ is a piecewise linear map from one interval $(a,b)$ to another, then the \emph{initial slope} of $k$ is the derivative of $k$ at any point between $a$ and the first node of $k$, and the \emph{final slope} of $k$ is the derivative of $k$ at any point between the final node of $k$ and $b$. The invariant $\Sigma_2$ is a list of positive real numbers. The entry for a bump domain $I_j=(p_{j-1},p_j)$ is the value of the initial slope of $f$ in $I_j$. 

Finally, the invariant $\Sigma_3$ is a list of equivalence classes of \emph{finite functions}. We calculate the entry for a bump domain $I_j=(p_{j-1}, p_j)$ as follows. Suppose first of all that $\Sigma_1=1$ in $I_j$. Define, for $x\in I_j$, the \emph{slope ratio} $f^*(x)=\frac{f'_+(x)}{f'_-(x)}$. Thus $f^*(x)=1$ except when $x$ is a node of $f$. Now define
\[
\phi_{f,j}: I_j\to \R, \ \ x\mapsto \prod_{n=-\infty}^\infty f^*(f^n(x)).
\]
Since $f$ has only finitely many nodes, only finitely many terms of this infinite product are distinct from $1$.
Let $p$ be the smallest node of $f$ in $I_j$ and let $p_\ast$ be the smallest node of $f$ in $I_j$ such that $\phi_{f,j}(p_\ast)\neq 1$ (such a node must exist). Define, for $s\in[0,1)$,
\[
\psi_{f,j}(s)=\phi_f(\lambda^s(r-p_{j-1})+p_{j-1}).
\]
Here $\lambda$ is the entry in $\Sigma_2$ corresponding to $I_j$ and $r$ is any point in the interval $(0,p)$ which satisfies the formula $r=f^{n}(p_\ast)$ for $n$ some negative integer. 

The function $\psi_{f,j}$ is a \emph{finite function}, that is, a function $[0,1)\to \mathbb{R}^+$ which takes the value $1$ at all but finitely many values. In our definition of $\psi_{f,j}$, we have chosen a value for $r$ which guarantees that $\psi_{f,j}(0)\neq 1$; we can do this by virtue of \cite[Lemma 4.4]{BrSq01}. 

The entry for $\Sigma_3$ corresponding to $I_j$ is the equivalence class $[\psi_{f,j}]$, where two finite functions $c_1$ and $c_2$ are considered equivalent if $c_1 = c_2\circ \rho$ where $\rho$ is a translation of $[0,1)$ modulo 1. If $f(x)<x$ for each $x\in I_j$ then the entry for $\Sigma_3$ corresponding to $I_j$ is the equivalence class $[\psi_{f^{-1}, j}]$.

\section{The definition of $\Delta$}\label{S: Delta}

The quantity $\Delta$ will also be a list, this time a list of equivalence classes of tuples of real numbers. To begin with we need the concept of a {\it minimum cornered function}.

\subsection{The minimum cornered function}
 
Take $f\in \text{PL}^+(0,1)$ and focus on the restriction of $f$ to one of its bump domains $D=(a,b)$. We adjust one of the definitions of Brin and Squier \cite{BrSq01}: for us, a \emph{cornered function} in $\text{PL}^+(0,1)$ is an element $l$ which has a single bump domain $(a,b)$ and which satisfies the following property: $\Sigma_1$ takes value $1$ (resp. $-1$) in relation to $(a,b)$ and there exists a point $x\in (a,b)$ such that all nodes of $l$ which lie in $(a,b)$ lie in $(x,l(x))$ (resp. $(l(x),x)$). We will sometimes abuse notation and consider such a cornered function as an element of $\text{PL}^+(a,b)$.

We say that a cornered function $l$ {\it corresponds to} a finite function $c$ if $\psi_l = c$ (we drop the subscript $i$ here, since there is only one bump domain). Roughly speaking this means that the first node of $l$ corresponds to $c(0)$. For a given initial slope $\lambda$ there is a unique cornered function in $\text{PL}^+(a,b)$ such that $\psi_l=c$ (this can be deduced from the proof of Lemma \ref{L: elementary}).

Now let $c:[0,1)\to \R$ be a finite function such that $[c]$ is the entry in $\Sigma_3$ associated with $D$. Within this equivalence class $[c]$ we can define a {\it minimum finite function} $c_m$ as follows. First define $C=\{c_1\in[c] | c_1(0)\neq 1\}$ and define an ordering on $C$ as follows. Let $c_1,c_2\in C$ and let $x$ be the smallest value such that $c_1(x)\neq c_2(x)$. Write $c_1<c_2$ provided $c_1(x)<c_2(x)$. We define $c_m$ to be the minimum function in $C$ under this ordering.

Suppose that $\lambda$ is the entry in $\Sigma_2$ associated with $D$. Suppose that $l$ is the cornered function in $\text{PL}^+(a, b)$, with initial slope $\lambda$, which corresponds to $c_m$. We say that $l$ is the {\it minimum cornered function} associated with $f$ over $D$. 

\subsection{The quantity $\Delta$}
Suppose now that  $f\in\ThF$. A \emph{bump chain} is a subsequence $I_t,I_{t+1},\dots,I_u$ of \eqref{E: intervals} such that each interval is a bump domain, and of the points $p_{t-1},p_t,\dots,p_u$ only $p_{t-1}$ and $p_{u}$ are dyadic. Thus $I_1, I_2,\dots, I_m$ can be partitioned into bump chains and open intervals of fixed points of $f$ (which have dyadic numbers as end-points). 

In \cite{BrSq01}, conjugating functions in $\text{PL}^+(0,1)$ are constructed by dealing with one bump domain at a time. We will construct conjugating functions in $\ThF$ by dealing with one bump \emph{chain} at a time. Consider a particular bump chain $D_1, \dots, D_s$ and let $f_j$ be the restriction of $f$ to $D_j=(a_j, b_j)$. 

According to \cite[Theorem 4.18]{BrSq01}, the centralizer of $f_j$ within $\text{PL}^+(a_j, b_j)$ is an infinite cyclic group generated by a root $\widehat{f_j}$ of $f_j$. We define $\lambda_j$ to be the initial slope of $\widehat{f_j}$ and $\mu_j$ to be the final slope of $\widehat{f_j}$. (Let $m_j$ be the integer such that  $\widehat{f_j}^{m_j}=f_j$; then $\lambda_j$ and $\mu_j$  are the positive $m_j$th roots of the initial and final slopes of $f_j$.)

Next, let $k_j$ be a member of $\text{PL}^+(a_j, b_j)$ that conjugates $f_j$ to the associated minimum cornered function, $l_j$, in $\text{PL}^+(a_j,b_j)$. Thus $k_j$ is some function satisfying $k_j f_j k_j^{-1}=l_j$. Let $\alpha_j$ be the initial slope of $k_j$ and let $\beta_j$ be the final slope. 

Consider the equivalence relation on $\mathbb{R}^s$ such that $(x_1,\dots,x_s)$ is equivalent to $(y_1,\dots,y_s)$ if and only if there are integers $m,n_1,\dots,n_s$ such that
\begin{align*}
2^mx_1 &=\lambda_1^{n_1}y_1\\
\mu_1^{n_1}x_2 &=\lambda_2^{n_2}y_2\\
\mu_2^{n_2}x_3 &=\lambda_3^{n_3}y_3\\
\vdots  &\\
\mu_{s-1}^{n_{s-1}}x_s &=\lambda_s^{n_s}y_s.
\end{align*}
It is possible to check whether two $s$-tuples of real numbers are equivalent according to the above relation in a finite amount of time because the $\lambda_i$ and $\mu_j$ are rational powers of $2$. We assign to the chain $D_1,\dots,D_s$ the equivalence class of the $s$-tuple
\[
\left(\frac{\alpha_1}{w_1},\frac{\alpha_2}{w_2}\frac{w_1}{\beta_1}, \dots , \frac{\alpha_s}{w_s}\frac{w_{s-1}}{\beta_{s-1}}\right)
\]
where $w_j=b_j-a_j$. We define $\Delta_f$ to consist of an ordered list of such equivalence classes; one per bump chain.

\section{Proof of Proposition~\ref{P: main}}\label{S: invariant}

We prove Proposition~\ref{P: main} after the following elementary lemma.
 
\begin{lemma}\label{L: inF}
Let $f$ and $g$ be maps in $\ThF$, and let $h$ be an element of $\textup{PL}^+(0,1)$ such that $hfh^{-1} = g$. Let $D=(a,b)$ be a bump domain of $f$ and suppose that the initial slope of $h$ in $D$ is an integer power of $2$. Then all slopes of $h$ in $D$ are integer powers of $2$ and all nodes of $h$ in $D$ occur in $\Z[\frac12]$.
\end{lemma}
\begin{proof}
Let $(a,a+\delta)$ be a small interval over which $h$ has constant slope; suppose that this slope is greater than $1$. We may assume that $f$ has initial slope greater than $1$ otherwise replace $f$ with $f^{-1}$ and $g$ with $g^{-1}$. Now observe that $hf^nh^{-1}=g^n$ for all integers $n$ and so $h = g^nhf^{-n}$.

Now, for any $x\in (a,b)$ there is an interval $(x,x+\epsilon)$ and an integer $n$ so that $f^{-n}(x,x+\epsilon)\subset (a,a+\delta)$. Then the equation $h=g^nhf^{-n}$ implies that, where defined, the derivative of $h$ over $(x,x+\epsilon)$ is an integral power of $2$. Furthermore any node of $h$ occuring in $(x,x+\epsilon)$ must lie in $\Z[\frac12]$ as required.

If $h$ does not have slope greater than $1$ then apply the same argument to $h^{-1}$ using the equation $h^{-1}gh=f$.
\end{proof}

We have two elements $f$ and $g$ of $\ThF$ and a third element $h$ of $\text{PL}^+(0,1)$ such that $hfh^{-1}=g$. We use the notation for $f$ described in the previous section, such as the quantities $I_j$, $p_j$, $f_j$, $\widehat{f_j}$, $k_j$, $l_j$, $\alpha_j$, $\beta_j$, $w_j$, $\lambda_j$, and $\mu_j$. We need exactly the same quantities for $g$, and we distinguish the quantities for $g$ from those for $f$ by adding a $'$ after each one. In particular, we choose a bump chain $D_1,\dots,D_s$ of $f$ and define $D'_i=h(D_i)$ for $i=1,\dots,s$. Note that $D'_1,\dots,D'_s$ are bump domains but need not form a bump chain for $g$ according to our assumptions, because $h$ is not necessarily a member of $\ThF$.

Let the function $h_i=h|_{D_i}$ have initial slope $\gamma_i$ and final slope $\delta_i$. Let $u$ be the member of $\textup{PL}^+(0,1)$ which, for $i=1,\dots,m$, is affine when restricted to $I_i$, and maps this interval onto $I'_i$. Notice that, restricted to ${D'}_i$, $u l_i u^{-1}$ is a cornered function which is conjugate to ${l'}_i$ (by the map $k'_ih_ik_i^{-1}u^{-1}$), and which satisfies $\psi_{{l'}_i}=\psi_{u l_i u^{-1}}$. Therefore $u l_i u^{-1} = {l'}_i$.  Combine this equation with the equations $k_if_ik_i^{-1}=l_i$, ${k'}_ig_i{k'}_i^{-1}={l'}_i$, and $h_if_ih_i^{-1}=g_i$ to yield
\[
(k_i^{-1}u^{-1} k'_ih_i)f_i (k_i^{-1}u^{-1} k'_ih_i)^{-1} = f_i.
\]
Therefore $k_i^{-1}u^{-1} k'_ih_i$ is in the centralizer of $f$, so there is an integer $N_i$ such that 
\[
h_i = (k'_i)^{-1}  u k_i {\widehat{f_i}}^{N_i} 
\]
for each $i=1,\dots,s$. Then by comparing initial and final slopes in this equation we see that 
\begin{equation}\label{E: gamma}
\gamma_i =\lambda_i^{N_i}\frac{\alpha_i}{w_i}\frac{w'_i}{\alpha'_i}, \quad 
\delta_i =\mu_i^{N_i}\frac{\beta_i}{w_i}\frac{w'_i}{\beta'_i}
\end{equation}
for $i=2,\dots,s$. We are now in a position to prove Proposition~\ref{P: main}.

\begin{proof}[Proof of Proposition~\ref{P: main}]
Suppose that $h\in \ThF$. Then there are integers $M_1,\dots,M_s$ such that $\gamma_1=2^{M_1}$ and $\gamma_i=\delta_{i-1}=2^{M_i}$ for $i=2,\dots,s$. Substituting these values into \eqref{E: gamma} we see that
\[
2^{M_1}\frac{\alpha'_1}{w'_1} = \lambda_1^{N_1}\frac{\alpha_1}{w_1},\quad \mu_{i-1}^{N_{i-1}}\frac{\alpha'_i}{w'_i}\frac{w'_{i-1}}{\beta'_{i-1}} = \lambda_i^{N_i}\frac{\alpha_i}{w_i}\frac{w_{i-1}}{\beta_{i-1}},
\]
for $i=2,\dots,s$, as required.

Conversely, suppose that $\Delta_f=\Delta_g$. We modify $h$ so that it is a member of $\ThF$. If $I_j$ is an interval of fixed points of $f$ then modify $h_j$ so that it is any  piecewise linear map from $I_j$ to $I'_j$ whose slopes are integer powers of $2$, and whose nodes occur in $\Z[\frac12]$. (It is straightforward to construct such maps, see  \cite[Lemma 4.2]{CaFlPa96}.) 

Now we modify $h$ on a bump chain $D_1,\dots, D_s$. Since $\Delta_f=\Delta_g$ we know that there are integers $m$ and $n_1,\dots,n_s$ such that, for $i=2,\dots,s$,
\begin{equation}\label{E: equiv}
2^m \frac{\alpha_1}{w_1} = \lambda_1^{n_1}\frac{\alpha'_1}{w'_1},\quad
\mu_{i-1}^{n_{i-1}}\frac{\alpha_i}{w_i}\frac{w_{i-1}}{\beta_{i-1}} = \lambda_i^{n_i}\frac{\alpha'_i}{w'_i}\frac{w'_{i-1}}{\beta'_{i-1}}.
\end{equation}

Consider the piecewise linear map $h'_i:D_i\rightarrow h_i(D_i)$ given by $h'_i=h_i\widehat{f_i}^{-n_i-N_i}$. The initial slope $\gamma'_i$ of $h'_i$ is $\gamma_i \lambda_i^{-n_i-N_i}$ and the final slope $\delta'_i=\delta_i \mu^{-n_i-N_i}$. From \eqref{E: gamma} and \eqref{E: equiv} we see that
\[
\gamma'_1=2^{-m},\quad \gamma'_i=\delta'_{i-1}.
\]
for $i=2,\dots,s$. We modify $h$ by replacing $h_i$ with $h'_i$ on $D_i$. Then $h$ does not have a node at any of the end-points of $D_1,\dots,D_s$ other than the first and last end-point. By Lemma~\ref{L: inF}, the nodes of $h_1$ occur in $\Z[\frac12]$ and the slopes of $h_1$ are all powers of $2$. Since the initial slope of $h_1$ coincides with the final slope of $h_1$, the same can be said of $h_2$. Similarly, for $i=2,\dots, s$, the initial slope of $h_i$ coincides with the final slope of $h_{i-1}$. We repeat these modifications for each bump chain of $f$; the resulting conjugating map is a member of $\ThF$.
\end{proof}

\section{Calculating $\alpha_i$ and $\beta_i$}\label{S:calc1}

It may appear that, in order to calculate $\Delta$, it is necessary to construct various conjugating functions. In particular to calculate $\alpha_i$ one might have to construct the function  in $\text{PL}^+(a_i, b_i)$ which conjugates $f_i$ to the conjugate minimum cornered function in $\text{PL}^+(a_i, b_i)$.

It turns out that this is not the case. The values for $\alpha_i$ and $\beta_i$ can be calculated simply by looking at the entries in $\Sigma_1, \Sigma_2$ and $\Sigma_3$ which correspond to $D_i$. In this section we give a formula for $\alpha_i$; we then observe how to use the formula for $\alpha_i$ to calculate $\beta_i$. 

In what follows we take $f$ to be a function in $\text{PL}^+(a, b)$ such that $f(x)\neq x$ for $x\in(a,b)$. Let $l$ be the minimum cornered function which is conjugate to $f$ in $\text{PL}^+(a, b)$.

\subsection{Calculating $\alpha_i$}
Suppose first that $f(x)>x$ for $x\in(a,b)$. Let $y_j$, for $j=0,\dots,t$ be the points at which the finite function $\psi_f$ does not take value $1$; let $\psi_f$ take the positive value $z_j$ at the point $y_j$ and assume that $0=y_0<y_1<\dots<y_t<1$. We will denote $\psi_f$ by $c_t$ and define $c_j=c_t(x+y_{j+1})$. Then $c_j$ is a translation of $c_t$ under which $y_j$ is mapped to the last point of $c_j$ which does not take value $1$.

Let $u_j$ be the cornered function corresponding to $c_j$ and let $x_j$ be the final node of $u_j$. Note that $u_j$ is conjugate to $f$ and, for $j$ equal to some integer $n$, $u_j$ equals $l$, the minimum cornered function. Define the {\it elementary function} $h_{x,r}$ to be  the function which is affine on $(0,x)$ and $(x,1)$ and which has slope ratio $r$ at $x$. We define $\zeta_j$ to be the initial slope of the elementary function $h_{x_j,z_j}$.

Let $p$ be the first node of $f$ and let $q$ be the first node of $u_t$.

\begin{lemma}\label{l:is comp}
There exists $k$ in $\textup{PL}^+(a, b)$ such that $kfk^{-1}=l$ and the initial slope of $k$ is
\[
(\zeta_t\zeta_{t-1}\dots\zeta_{n+1})\left(\frac{q-a}{p-a}\right).
\]
\end{lemma}

Note that, in the formula just given, $p$ and $q$ stand for the $x$-coordinates of the corresponding nodes. Before we prove Lemma \ref{l:is comp} we observe that we can calculate values for the $\zeta_j$ and $q$ simply by looking at $\Sigma_2$ and $\Sigma_3$ and using the following lemma:

\begin{lemma}\label{L: elementary}
Let $l$ be a cornered function in $\textup{PL}^+(a,b)$ with initial slope $\lambda>1$, and suppose that the corresponding finite function $c$ takes the value $1$ at all points in $[0,1)$ except $0=s_0<s_1<\dots<s_k<1$, at which $c(s_i)=z_i$. Then the first node $q_0$ of $l$ is given by the formula
\begin{equation}\label{E: q_0}
q_0=a+\tfrac{(b-a)(1-[\lambda z_0\dotsb z_k])}{[\lambda(1-z_0)]+[\lambda^{s_1+1}z_0(1-z_1)]+\dots+[\lambda^{s_{k-1}+1}z_0\dotsb z_{k-2}(1-z_{k-1})]+[\lambda^{s_k+1}z_0\dotsb z_{k-1}(1-z_k)]},
\end{equation}
and the initial slope $\zeta$ of the elementary function $h_{q_k,z_k}$, where $q_k$ is the final node of $l$, is given by 
\begin{equation}\label{E: zeta}
\zeta = \frac{b-a}{\lambda^{s_k}(q_0-a)(1-z_k)+(b-a)z_k}.
\end{equation}
\end{lemma}
\begin{proof}
If $q_0,\dots,q_k$ are the nodes of $l$ we have equations
\begin{equation}\label{E: q_i}
\lambda^{s_i}(q_0-a)+a=q_i, \quad i=0,1,2,\dots,k.
\end{equation}
Define $q_{k+1}=b$ and let $\lambda_i$ be the slope of $l$ between the nodes $q_{i-1}$ and $q_{i}$ for $i=1,\dots,k+1$. Then $z_i=\lambda_{i}/\lambda_{i-1}$ for $i>1$, and we obtain
\begin{equation}\label{E: lambda_i}
\lambda_i = \lambda z_0\dots z_{i-1},\quad i=1,\dots,k+1.
\end{equation}
If we substitute \eqref{E: q_i} and \eqref{E: lambda_i} into the equation
\[
b-a = \lambda(q_0-a) + \lambda_1(q_1-q_0) + \lambda_2(q_2-q_1) + \dots + \lambda_{k+1}(b-q_k),
\]
then we obtain \eqref{E: q_0}. To obtain \eqref{E: zeta}, notice that $z_k\zeta$ is the final slope of $h_{q_k,z_k}$, therefore $b-a = \zeta(q_k-a) +z_k\zeta (b-q_k)$. Substitute the value of $q_k$ from \eqref{E: q_i} into this equation to obtain \eqref{E: zeta}.
\end{proof}

Before we prove Lemma \ref{l:is comp}, we make the following observation. Let $g$ be a function such that $g(x)>x$ for all $x\in(a,b)$ and suppose that $g$ has nodes $p_1<\dots< p_s$. Now let $h = h_{p_s, g^*(p_s)}$. Then $hgh^{-1}$ has nodes
$h(p_1),\dots, h(p_{s-1}), hg^{-1}(p_s)$ with $(hgh^{-1})^*$ taking on values $g^*(p_1), \dots g^*(p_t)$ at the respective nodes. If $ hg^{-1}(p_s)=h(p_i)$ for some $i$, then $(hgh^{-1})^*$ has value $g^*(p_i)g^*(p_s)$.

\begin{proof}[Proof of Lemma \ref{l:is comp}]
The formula given in Lemma \ref{l:is comp} arises as follows. We start by finding the conjugator from $f$ to the cornered function $u_t$; then we cycle through the cornered functions $u_j$ until we get to $u_n=l$. Thus the $\frac{q-a}{p-a}$ part of the formula arises from the initial conjugation to a cornered function, and the $\zeta_j$'s arise from the cycling.

Consider this cycling part first and use our observation above on the cornered functions, $u_j$: we have $h_{x_j, z_j}u_j(h_{x_j, z_j})^{-1}=u_{j-1}$. Thus in order to move from $u_t$ to $u_n$ we repeatedly conjugate by elementary functions with initial gradient $\zeta_t,\dots,\zeta_{n+1}$.

We must now explain why we can use $\frac{q-a}{p-a}$ for the first conjugation which moves from $f$ to $u_t$. It is sufficient to find a function which conjugates $f$ to $u_t$ and which is linear on $[a,p]$.

Consider the effect of applying an elementary conjugation to a function $f$ that is not a cornered function. Suppose that $f$ has nodes $p_1<\dots<p_s$. So $p=p_1$. We consider the effect of conjugation by an elementary function $h = h_{p_s, f^*(p_s)}$ as above. To reiterate, we obtain a function with nodes
$$h(p_1), ...., h(p_{s-1}), hf^{-1}(p_s)$$
Now observe that, since $f^*(p_s)<1$, $h(x)>x$ for all $x$ and $h$ is linear on $[a,p_s]$. So clearly $h$ is linear on the required interval. There are three possibilities:
\begin{itemize}
\item If $hf^{-1}(p_s) < h(p)$ then $f$ was already a cornered function; in fact $f=u_t$. We are done.
\item If $hf^{-1}(p_s) > h(p)$ then we simply iterate. We replace $f$ with $hfh^{-1}$, $p$ with $h(p)$ etc. We conjugate by another elementary function exactly as before. It is clear that the next elementary conjugation will be linear on $[a, h(p)]$ which is sufficient to ensure that the composition is linear on $[a,p]$.
\item If $hf^{-1}(p_s)=h(p)$ then we need to check if $hfh^{-1}$ is a cornered function. If so then $hfh^{-1}=u_t$, the corner function we require. If $hfh^{-1}$ is not a cornered function then we iterate as above, replacing $f$ with $hfh^{-1}$. It is possible that $h(p)$ will no longer be the first node of $hfh^{-1}$, but in this case we replace $p$ by $h(p_2)$. Since $[a,h(p_2]\supset [a,h(p)]$ this is sufficient to ensure that the composition is linear on $[a,p]$.
\end{itemize}

We can proceed like this until the process terminates at a cornered function. Since conjugating a non-cornered function by $h$ preserves $\psi$ we can be sure that we will terminate at $u_t$ as required. What is more the composition of these elementary functions is linear on $[a,p]$.
\end{proof}

Suppose next that $f(x)<x$ for all $x\in(a,b)$  and $kfk^{-1}=l$, a minimum cornered function. Observe that $kf^{-1}k^{-1} = l^{-1}$ and $f^{-1}(x)>x$ for all $x\in(a,b)$. We can now apply the formula in Lemma \ref{l:is comp}, replacing $f$ with $f^{-1}$ and $l$ with $l^{-1}$, to get a value for the initial slope of $k$. 

\subsection{Calculating $\beta_i$}
The method we have used to calculate $\alpha_i$ can also be used to calculate $\beta_i$. Define
$$\tau:[a,b]\to[a,b], x\mapsto b+a-x.$$
Now $\tau$ is an automorphism of $\text{PL}^+(a,b)$; the graph of a function, when conjugated by $\tau,$ is rotated $180^\circ$ about the point $(\frac{b+a}2, \frac{b+a}2)$. Consider the function $\tau f \tau$ and let $k$ be the conjugating function from earlier, so that $kfk^{-1}=l$. Then
$$(\tau k \tau) (\tau f \tau) (\tau k \tau)^{-1} = (\tau l \tau).$$

The initial slope of $\tau k \tau$ equals the final slope of $k$. Thus we can use the method outlined above -- replacing $f$ with $\tau f \tau$ and $l$ with $\tau l \tau$ -- to calculate the initial slope of $\tau k \tau$. Note that, for this to yield $\beta_i$, we must make an adjustment to the integer $n$ in the formula in Lemma \ref{l:is comp}: the function $\tau l \tau$ is not necessarily the {\it minimum} cornered function which is conjugate to $\tau f \tau$. Thus we choose $n$ to ensure that $l$ is minimum rather than $\tau l \tau$.

\section{Calculating $\lambda_i$ and $\mu_i$}\label{S:calc2}

Let $f$ be a fixed-point free element of $\text{PL}^+(a,b)$. Let $\widehat{f}$ be a generator of the centralizer of $f$ within $\text{PL}^+(a,b)$. The formula for $\Delta$ requires that we calculate the initial slope and the final slope of $\widehat{f}$. It turns out that this is easy---thanks to the work of Brin and Squier \cite{BrSq01}.

Let $c, c':[0,1)\to\R$ be finite functions. We say that $c'$ is the $p$-th root of $c$ provided that, for all $x\in [0,1)$, we have $c(x)=c'(px)$. The property of having a $p$-th root is preserved by the equivalence used to define $\Sigma_3$. Thus we may talk about the equivalence class $[c]$ having a $p$-th root, provided any representative of $[c]$ has a $p$-th root.

Now \cite[Theorem 4.15]{BrSq01} asserts that $f$ has a $p$-th root in $\text{PL}^+(a,b)$, for $p$ a positive integer, if and only if the single equivalence class in $\Sigma_3$ is a $p$-th power (following Brin and Squier we say that this class has \emph{$p$-fold symmetry}). What is more \cite[Theorem 4.18]{BrSq01} asserts that $\widehat{f}$ must be a root of $f$.

Thus if $p$ is the largest integer for which the single class in $\Sigma_3$ has $p$-fold symmetry then  $\widehat{f}$ is the $p$-th root of $f$. The initial slope of $\widehat{f}$ is the positive $p$-th root of the initial slope of $f$, and the final slope of $\widehat{f}$ is the positive $p$-th root of the final slope of $f$.

\section{Roots and centralizers}\label{S:rc}


We describe the roots and centralizers in $\ThF$ by extending the work of Brin and Squier on $\text{PL}^+(0,1)$. Let $f$ be a non-trivial element of $\ThF$ with bump domains $E_1, \dots, E_k$.  Let $f_i$ denote the restriction of $f$ to $E_i=(a_i, b_i)$ and, for $i=1,\dots,k$, define $m_i$ to be the integer such that the initial slope of $f_i$ is $2^{m_i}$. Define $p_i$ to be the largest integer such that the entry in $\Sigma_3$ for $E_i$ has $p_i$-fold symmetry.

We are interested in giving conditions for an element to be a root, or a centralizer, of $f$. We also want to give the structure of the following groups:
$$R(f)=\langle g\in \ThF: g^a = f \text{ for some } a\in\mathbb{Z}\rangle, \quad 
C(f)=\{g\in \ThF :gf=fg\}.$$

\begin{theorem}\label{T: root}
Let $f$ be a non-trivial element of $\ThF$. Given an integer $p$, there is an element $g$ of $\ThF$ such that $g^p=f$ if and only if $p$ divides each of $p_1, \dots, p_k, m_1, \dots, m_k$. Then $R(f)\cong\mathbb{Z}$.
\end{theorem}
\begin{proof}
Suppose $F$ contains an element $g$ such that $g^p=f$. Then $g$ shares the same bump domains as $f$ and \cite[Theorem 4.15]{BrSq01} implies that each entry of $\Sigma_3$ must have $p$-fold symmetry. Hence $p$ divides $p_i$ for $i=1,\dots, k$. Furthermore the initial slope $2^{n_i}$ of $g|_{E_i}$ satisfies $2^{n_ip}=2^{m_i}$. Therefore $p|m_i$. 

Conversely suppose that $p$ is an integer dividing each of $p_1, \dots, p_k, m_1, \dots, m_k$. Then \cite[Theorem 4.15]{BrSq01} implies that there exists a map $g$, which is a $p$-th root of $f$ in $\text{PL}^+(0,1)$. Again, $g$  shares the same bump domains as $f$, and on $E_i$, the initial slope $\gamma_i$ of $g|_{E_i}$ satisfies $\gamma_i^p=2^{m_i}$. Therefore $\gamma_i=2^{m_i/p}$.  Lemma~\ref{L: inF} applies to show that, within bump domains, all slopes of $g$ are powers of $2$, with all nodes in $D$ occurring in $\Z[\frac12]$. What is more, if $E_i$ and $E_{i+1}$ occur in the same bump chain in $f$ then the initial slope of $f$ in $E_{i+1}$ must equal the final slope of $f$ in $E_i$. Clearly this property will also transfer to $g$; hence $g$ is an element of $\ThF$. 

Now \cite[Theorem 4.15]{BrSq01} asserts that a $p$-th root of $f$ is unique in $\text{PL}^+(0,1)$. Hence $R(f)$ is generated by $g$ where $g$ is the smallest root of $f$ in $\ThF$; thus $R(f)\cong\mathbb{Z}$.
\end{proof}

Brin and Squier also proved in \cite[Theorem 4.18]{BrSq01} that the only maps in $\text{PL}^+(a,b)$ that commute with a fixed point free member $g$ of $\text{PL}^+(a,b)$ are roots of $g$.  We can use this to describe the maps in $\ThF$ which commute with $f$. 

For $a,b\in \mathbb{Z}[\frac12]\cap[0,1]$, let $C=(a,b)$ and define $\ThF_C$ to be the group consisting of those elements in $\ThF$ which fix $x$ for $x\not\in C$. Then \cite[Lemma 4.4]{CaFlPa96} implies that $\ThF_C \cong \ThF$ provided $b-a$ is a power of $2$. In fact it is easy to extend the proof of \cite{CaFlPa96} to prove that, even without this proviso, $\ThF_C \cong \ThF$: break $C$ into intervals $(p_{i-1}, p_i)$ for $i=1,\dots, s$ such that $p_0=a, p_s = b$ and $p_i - p_{i-1}$ is a power of $2$. Similarly, break $(0,1)$ into intervals $(q_{i-1}, q_i)$ for $i=1,\dots, s$ such that $q_0=0, q_s=1$ and $q_i-q_{i-1}$ is a power of $2$. Note that $p_i$ and $q_i$ are dyadic rational numbers for $i=0,\dots, s$. Now define the map $k:(0,1)\to(a,b)$ such that $k(q_i)=p_i$ for $i=1,\dots, s$ and $k$ is affine on the interval $(q_{i-1}, q_i)$. Then we define a map $\phi:\ThF\to\ThF_C$ as follows: for $f\in\ThF$, 
$$(\phi(f))(x) = 
\left\{\begin{array}{ll}
(kfk^{-1})(x), & x\in C; \\
x, & x\not\in C. 
\end{array}
\right.$$ 
It is easy to check that $\phi$ is an isomorphism from $\ThF$ to $\ThF_C$.

\begin{theorem}
Suppose that $g\in \ThF$ commutes with $f\in\ThF$. Then
\begin{itemize}
\item the restriction of $g$ to a bump chain of $f$ is the root of a restriction of $f^p$ for some integer $p$;
\item the restriction of $g$ to a maximal connected open set, $C$, of $\textup{fix}(f)$ is any element of $\ThF_C$.
\end{itemize}
Then $C_{\ThF}(f)\cong \ThF^a\times \mathbb{Z}^b$ where $a$ is the number of non-empty maximal connected open sets in $\textup{fix}(f)$ and $b$ is the number of bump chains of $f$.
\end{theorem}
\begin{proof}
Let $D_1, \dots, D_s$ be a bump chain of $f$. Now, by \cite[Theorem 4.18]{BrSq01}, we know that $g|_{D_i}$ is a root of $f|_{D_i}$ for $i=1,\dots, s$. In order for $g$ to lie in $\ThF$, the final slope of $g$ in $D_{i-1}$ must equal the initial slope in $D_i$ for $i=1,\dots, s-1$. Since the same is true of $f$ we conclude that there exist integers $p$ and $q$ such that $g|_{D_i}^q=f|_{D_i}^p$ for $i=1,\dots, s.$

Clearly if $C$ is a maximal connected open set in $\text{fix}(f)$ then $f|_C$ may coincide with any element of $\ThF_C$. The structure for $C_{\ThF}(f)$ follows easily.
\end{proof}

\bibliographystyle{plain}
\bibliography{reversible}
\end{document}